\documentclass[a4paper]{amsart}
\usepackage[
left=38truemm,
right=38truemm
]{geometry}
\usepackage{amsfonts, amssymb, amsmath, verbatim, amsthm, mathrsfs, amscd, enumerate, ascmac, fancyhdr, caption, hyperref, framed, subfigure}
\usepackage{multirow}
\usepackage{bbm}
\usepackage{tikz}
\usetikzlibrary{positioning,intersections, calc, arrows,decorations.markings, angles}
\colorlet{Green}{green!70!black!}
\theoremstyle{plain}
\newtheorem{theorem}{Theorem}
\newtheorem{lemma}[theorem]{Lemma}
\newtheorem*{lemma*}{Lemma}
\newtheorem{proposition}[theorem]{Proposition}
\newtheorem*{proposition*}{Proposition}

\theoremstyle{definition}

\newtheorem*{claim*}{Claim}
\theoremstyle{remark}

\def\d{\mathrm{d}}
\def\ae{\textrm{a.e.}}

\def\C{\mathfrak C}

\def\D{\mathfrak{D}}

\author{Chu-hee Cho}
\address[Chu-hee Cho] {Department of Mathematical Sciences and RIM, Seoul National University, Seoul 08826, Republic of Korea}
\email{akilus@snu.ac.kr}

\author{Shobu Shiraki}
\address[Shobu Shiraki] {Graduate School of Mathematics, Nagoya University, Furocho, Chikusa-ku, Nagoya 464-8602, Japan}
\email{shiraki.shobu@math.nagoya-u.ac.jp}

\begin{document}
\title[
]
{
Dimension of divergence sets of oscillatory integrals with concave phase
}
\begin{abstract}
We study the Hausdorff dimension of the sets on which the pointwise convergence of the solutions to the fractional Schr\"odinger equation $e^{it(-\Delta)^\frac m2}f$ fails when $m\in(0,1)$ in one spatial dimension. The poinwise convergence along a non-tangential curve and a set of lines are also considered, where we find a different nature compared to the case when $m\in(1,\infty)$. 

\end{abstract}

\maketitle
\section{Introduction}

For $n\in\mathbb N$ and $m\in(0,\infty)$, the solution to the Schr\"odinger-type equation on $\mathbb R^n\times\mathbb R$
\[
\begin{cases}
i\partial_tu(x,t)+(-\Delta)^\frac m2 u(x,t)=0,\\
u(x,0)=f(x)
\end{cases}
\]
is (formally) given by
\[
u(x,t)
=
e^{it(-\Delta)^{\frac m2}}f(x)
=
(2\pi)^{-n}\int_{\mathbb R^n}e^{i(x\cdot\xi+t|\xi|^m)}\widehat{f}(\xi)\,\d\xi.
\]
Let $m=2$ in which case, it is the classical and standard Schr\"odinger equation in quantum mechanics. A fundamental problem is to determine the least smoothness $s$ of the initial data for which the pointwise convergence to the solutions 
\begin{equation}\label{pw}
\lim_{t\to 0} e^{it(-\Delta)^\frac m2}f(x)=f(x)\quad \ae
\end{equation}
is guaranteed, where the initial data $f$ is taken from a (inhomogeneous) Sobolev space $H^s(\mathbb R^n)$ defined by
\[
\|f\|_{H^s(\mathbb R^n)}
=
\|(1-\Delta)^\frac s2 f\|_{L^2(\mathbb R^n)}.
\]
This problem is often called Carleson's problem and traced back to work by Carlson \cite{Cr80} and Dahlberg--Kenig \cite{DK82} for the Schr\"odinger equation in the early 1980s. They solved it completely in one spacial dimension: The pointwise convergence \eqref{pw} with $(n,m)=(1,2)$ holds if and only if $s\geq\frac14$. The higher dimensional cases are more difficult to study. After many authors' contributions (see \cite{Sj87, Vg88, Br95, Lee06, Br13, LR17, LR19a}, for example), Bourgain \cite{Br16} (for which a nice expository paper by Pierce \cite{Pr20} is helpful), Du--Guth--Li \cite{DGL17}, and Du--Zhang \cite{DZ18} finally proved that $s=\frac12-\frac{1}{2(n+1)}$ is the critical regularity, although it is still unknown whether the pointwise convergence is true exactly at the critical point.

The fractional Schr\"odinger equations ($m\in(0,\infty)\backslash\{1\}$) are natural generalizations not only in mathematics but also in physics. For instance, L\'evy stochastic process generalizes the Gaussian process or Wiener stochastic process, in which context the parameter $m$ may be called L\'evy index. 
Furthermore, the fractional Schr\"odinger equations have connections in optics related to Airy beam as well as water wave equations with surface tension. The reader may visit \cite{Ls00, Ls02, GH11, IP14, Longhi15} and references therein.
In view of the pointwise convergence problem, the case of the fractional Schr\"odinger equations are of interest in their own right.

In the case when $m > 1$, the regularity appears to be the same as that for the standard Schr\"odinger equation. In one dimension, the regularity is known to be independent of $m$ \cite{Sj87, KPV91}. In higher dimensions, it is known that $s > \frac{1}{2} - \frac{1}{2(n-1)}$ is at least sufficient \cite{CK18}, but the necessary condition remains an open problem.


Sj\"ogren--Sj\"olin \cite{SS89} and Barcel\`o--Bennett--Carbery--Rogers \cite{BBCR11} introduced a refinement of the pointwise convergence problem; measuring the Hausdorff dimensions of the so-called divergence sets. The divergence set $\mathfrak{D}(f)$ for each $f\in H^s(\mathbb R^n)$ is the set on which the pointwise convergence fails, namely, 
\[
\mathfrak{D}(f)
:=
\{x\in\mathbb R^n:\lim_{t\to0}e^{it(-\Delta)^\frac m2}f(x)\not=f(x)\}.
\]
Of course, quantifying the Hausdorff dimension of the divergence set is only meaningful to the smooth regularity $s$ for which the poinwise convergence holds; otherwise, the dimension is trivially full. 
In one dimension, Barcel\`o--Bennett--Carbery--Rogers \cite{BBCR11} revealed that 
\[
\sup_{f\in H^s(\mathbb R)}\dim_H\mathfrak D(f)=1-2s
\]
for $s\in(\frac14,\frac12]$ (For the lower bound they used the results for the Bessel potential due to \v{Z}ubrin\'{c} \cite{Zb02}). One may note that there is an interesting jump at $s=\frac14$. For higher dimensions, there are some partial results but many cases are still wide open. The interested readers may visit \cite{Mt95,BBCR11,BR12,LR17,LR19a,DGLZ18,EP22b}. 

When $m\in(0,1]$, it seems fairly different in nature, and less well known. 
The case when $m = 1$ corresponds to the wave equation. It has been discovered that $e^{it\sqrt{-\Delta}}f$ converges to $f$ almost everywhere for all $f \in H^s(\mathbb{R}^n)$ if $s > \frac{1}{2}$, as shown by Cowling \cite{Cw82}, and it fails otherwise, as demonstrated by Walther \cite{Wl97} (see also \cite{RV08}). Several studies have conducted deeper analyses regarding its divergence set. For instance, see \cite{BBCR11, LR19b, HKL21}. While the size of the divergence sets has been understood well in lower dimensions ($n = 1, 2, 3$), the problem in general remains open.

When $m\in(0,1)$, Cowling also noted in \cite{Cw82} that \eqref{pw} holds at least if $s>\frac m2$ in general $n \geq1$ and later Walther \cite{Wl95} proved that, in one dimension, \eqref{pw} holds if $s>\frac m4$ by considering the corresponding maximal inequality whose failure is also shown if $s < \frac m4$. The next question may be how big the size of divergence sets can be and our first result gives a reasonable upper bound in one dimension.


\begin{theorem}\label{t:concave}
Let $n=1$ and $m\in(0,1)$. Then 
\[
\sup_{f\in H^s(\mathbb R)}\dim_H\mathfrak D(f)
\leq
\max\left\{1-2s,\, \frac12+\frac{1-4s}{2(1-m)}\right\}
\]
whenever $s\in(\frac m4,\frac12)$.
\end{theorem}

The $\alpha$-dimensional measure $\mu$ and the corresponding maximal estimate associated with $\mu$ are central to proving Theorem \ref{t:concave}. The $\alpha$-dimensional measure is given by 
\[
\sup_{x\in\mathbb R^n,r>0}\frac{\mu(B(x,r))}{r^\alpha}
<
\infty.
\]
An $n$-dimensional ball of radius $r$ centered at $x$ is denoted by $B(x,r)$. We collect all $\alpha$-dimensional measures supported in the unit ball and denote such collection by $\mathcal M^\alpha$. 
By the standard argument (see \cite{CS21} for example), Theorem \ref{t:concave} follows from the following: for $q\geq2$
\begin{equation}\label{e:max}
    \|e^{it(-\Delta)^{\frac m2}}f\|_{L_x^q(\mathbb I,\d\mu)L_t^\infty(\mathbb I)}
    \lesssim
    \|f\|_{H^s}
\end{equation}
holds for all $\mu\in\mathcal M^\alpha$ and $f\in H^s$ and whenever
\begin{equation}\label{s:verline}
s
>
\max\left\{
\frac12-\frac m4-\frac{(1-m)\alpha}{q},
\frac12-\frac{\alpha}{q}
\right\}.
\end{equation}
The result is sharp in the sense that one can find an initial data and $\alpha$-dimensional measure such that \eqref{e:max} fails if $s<
\max\left\{
\frac12-\frac m4-\frac{(1-m)\alpha}{q},
\frac12-\frac{\alpha}{q}
\right\}$.
We discuss this in Section \ref{s:vertical}. 
The proof is constructed in the spirit of \cite{Sh19} where the second author dealt with a similar situation for $m>1$. Overall, the case with $m\in (0,1)$ seems much more delicate on which we make further comments after the proof. 
It is worth noting that the dimension of the divergence sets now continuously varies (Figure \ref{fig:dim}).

\begin{figure}[t]
\begin{center}
\begin{tikzpicture}[scale=1.5]

\draw [->](0,0)--(4.5,0) node [below] {$s$};
\draw [->](0,0)--(0,3.5) node [left] {$\dim_H\mathfrak D$};

\node [below left] at (0,0) {$O$};
\draw (1.4,3)--(0,3);
\draw [cyan] (4,0)--(2,1.5);
\draw [magenta] (2,1.5)--(1.4,3);
\draw [dotted] (1.4,3)--(1.4,0) node [below] {$\frac m4$};
\draw [dotted] (2,3)--(2,0) node [below] {$\frac14$};
\node [below] at (4,0) {$\frac12$};
\draw [dotted] (2,1.5)--(0,1.5) node [left] {$\frac12$};
\node [left] at (0,3) {$1$};

\draw [dashed](1.4,3)--(2,3);
\draw [dashed,cyan] (2,1.5)--(0,3);
\draw [dashed,magenta] (2,1.5)--(2.6,0);
\node [below] at (2.6,0) {$\frac12-\frac m4$};

\end{tikzpicture}
\caption{The upper bound of the divergence sets when $m\in(0,1)$.} \label{fig:dim}
\end{center}
\end{figure}

More variations regarding the pointwise convergence problem have been studied. 
Lee--Rogers \cite{LR12} and Lee--Vargas with the first author \cite{CLV12} considered pointwise convergence of solutions to the (standard) Schr\"odinger equations along a curve $\gamma$. When $n=1$, typically, it is given by
\begin{equation}\label{curve}
\gamma(x,t)=x-\theta t^\kappa, \qquad (x,t)\in\mathbb R\times\mathbb R, \quad \theta>0
\end{equation}
for $\kappa>0$. We say $\gamma(x,t)$ is \textit{non-tangential} if $\kappa \in [1,\infty)$ and \textit{tangential} if $\kappa\in (0,1)$. Figure \ref{fig:curves} illustrates the various paths that are discussed here.
The motivation in \cite{LR12} was to understand the pointwise convergence properties of the Shcr\"odinger equation with the harmonic oscillator. It turns out that it is equivalent to the one for the standard Schr\"odinger equation along a non-tangential curve $\sqrt{1+t^2}x$, which is classified the same as \eqref{curve}. In \cite{LR12} it was shown that the pointwise convergence of the Schr\"odinger equation\footnote{At least for a curve given by \eqref{curve}, our method generalizes the result to the fractional Schr\"odinger setting when $m>1$.} along a non-tangential curve given by \eqref{curve} (more generally $\gamma\in C^1$) 
\begin{equation}\label{pw along curve}
\lim_{\substack{(y,t)\to(x,0)\\ y=\gamma(x,t)}} e^{it(-\Delta)^{\frac m2}} f(y)
=
f(x)\quad \ae
\end{equation}
holds for all $f\in H^s(\mathbb R)$ if $s\geq\frac14$. This result reflects our intuition since a non-tangential curve ``looks like"  a vertical line around $t=0$ where pointwise convergence matters. (This can be more relaxed. For instance, there is no difference between a vertical line and a tilted line \eqref{curve} with $\kappa=1$ as long as $m>1$.) A vertical line, indeed, can be formally seen as \eqref{curve} with $\kappa=\infty$.

When the curve is tangential, it can be different from a vertical line. In fact, Lee, Vargas, and the authors \cite{CLV12,CS21} showed that \eqref{pw along curve} holds for all $f\in H^s(\mathbb R)$ provided that the curve is given by \eqref{curve} and $s>\max\{\frac14,\frac{1-m\kappa}{2}\}$.
The sizes of the divergence sets have already been estimated as 
\[
\dim_H \mathfrak D(f\circ\gamma)\leq \max\left\{1-2s,\frac{1-2s}{m\kappa}\right\}
\]
for $m>1$, $\kappa\in(0,1)$, and $s\in(\frac14,\frac12)$ (see \cite{CL14,CS21}). 
There are some partial results in higher dimensions \cite{LR12, LW18} but it is widely open at this moment. 

\begin{figure}[t]
\begin{center}
\begin{tikzpicture}[scale=0.8]

\tikzset{->-/.style={decoration={
  markings,
  mark=at position .5 with {\arrow{stealth}}},postaction={decorate}}}


 \draw [->](-1,0)--(6,0) node [below] {$y$};
 \draw [->](0,-0.5)--(0,3.5) node [left] {$t$};

 \draw [color=blue,->-] (5.5,3) to [out =-100, in=5] (2.5,0) ;
 \draw [color=orange,->-] (2.5,3)--(2.5,0);
 \draw [color=Green, ->-] (4, 3)--(2.5,0); 
 \draw [color=red, ->-] (-1,3) to [out =-10, in=95] (2.5,0);
  
 \node [below]at (2.5,0){$x$};

\end{tikzpicture}
\caption{Typical four kinds of paths; \textcolor{red}{non-tangential curve}, \textcolor{orange}{vertical line}, \textcolor{Green}{tilted (angled) line}, \textcolor{blue}{tangential curve}.} \label{fig:curves}
\end{center}
\end{figure}

When $m=1$, there is no space to make it interesting in one spacial dimension since the Hausdorff dimension of the divergence set is $0$ for $s>\frac n2$. On the other hand, the case $0<m<1$ seems surprisingly sensitive to the directions of the paths of convergence. Recently, Yuan--Zhao \cite{YZ21} showed that the pointwise convergence along a curve $\gamma$ given by \eqref{curve} with $0<\kappa\leq1$ holds if $s>\max\{\frac12-\frac m4,\frac{1-m\kappa}{2}\}$. They also showed that the Hausdorff dimension of the divergence sets is bounded above by $\frac{1-2s}{m\kappa}$ whenever $s\in(\frac12-\frac m4, \frac12)$.
When $\kappa=1$, the curve $\gamma$ becomes a (tilted) line with its angle (from the vertical line) $\theta$ indicating how different from a vertical line. By comparison with Walther's result introduced earlier, their results teach us that the path along a tilted line with a small angle is somehow very different from the vertical line in this context. The following theorem reveals that this phenomenon can be even severer.

\begin{theorem}\label{t:nontan}
Let $n=1$ and $m\in (0,1)$, $\kappa\geq 1$ and $\gamma$ given by \eqref{curve}. Then, the pointwise convergence along the curve $\gamma$ holds that \eqref{pw along curve} for all $f\in H^s(\mathbb R)$ if $s>\frac12-\frac m4$. Furthermore, we have
\[
\dim_H\mathfrak{D}(f\circ\gamma)
\leq
\frac{1-2s}{m}
\]
whenever $s\in(\frac12-\frac m4,\frac12)$.
\end{theorem}

It is possible to replace the curve $\gamma$ by the one in more general class (Remarks in Section \ref{s:non-tangential}).
Theorem \ref{t:nontan} is again obtained by the maximal(-in-time) inequality with respect to $\d\mu$. In section \ref{s:non-tangential}, we show that for $q\geq2$,
\begin{equation}\label{e:max general}
\|e^{it(-\Delta)^{\frac m2}} f\circ\gamma \|_{L_x^q(\mathbb I,\d\mu)L_t^\infty(\mathbb I)}
\leq
C\|f\|_{H^2(\mathbb R)}
\end{equation}
holds for all $\mu\in\mathcal M^\alpha$ and $f\in H^s(\mathbb R)$ if $s>\max\{\frac12-\frac m4,\frac12-\frac{m\alpha}{q}\}$. This is sharp in the sense that there exist an initial data and $\d\mu(x)$ such that \eqref{e:max general} fails if $s<\max\{\frac12-\frac m4,\frac12-\frac{m\alpha}{q}\}$. In particular, the special case when $\d\mu(x)=\d x$ (implying $\alpha=1$) directly gives us \eqref{pw along curve}. As \eqref{e:max general} shows, the property of pointwise convergence along a curve given by \eqref{curve} with a large $\kappa$ suddenly turns nicer once the curve becomes its superposition as $\kappa\to \infty$ (i.e. a vertical line). This interesting shift is even significant when $m$ is close to $0$. 
In fact, one aspect to understand why a non-tangential curve differs so much from a vertical line when $m\in(0,1)$ is that, formally speaking, if we consider $m=0$ and $\kappa>0$, we have 
\[
\|\sup_{t\in \mathbb I}|e^{it(-\Delta)^{0}}f(x)\|_{L^2(0,1)}
=
\|f\|_{H^0}
\]
on one hand, and by a change of variables
\[
\|\sup_{t\in \mathbb I}|e^{it(-\Delta)^{0}} f(x-\theta t^\kappa)|\|_{L^2(0,1)}
\lesssim
\|\sup_{t\in \mathbb R}|e^{it(-\Delta)^{\frac 12}} f(x)|\|_{L^2(0,1)}
\lesssim
\|f\|_{H^{\frac12+\varepsilon}}
\]
for arbitrary small $\varepsilon>0$ on the other hand. The cases when $m\in(0,1)$ capture a similar phenomenon of this shift in milder manners.

Theorem \ref{t:nontan} completes the study (except on the critical point) by providing reasonable sufficient regularities that guarantee pointwise convergence along a typical curve given by \eqref{curve} in one spatial dimension. Here is a brief summary: 
When $m>1$, the pointwise convergence \eqref{pw along curve} holds for all $f\in H^s(\mathbb R)$ if 
\[
s>\max\Big\{\frac14,\frac{1-m\kappa}{2}\Big\}\quad \text{for $\kappa>0$}.
\]
Furthermore, for such $s$, its Hausdorff dimension is bounded above by
\[
\dim_H\D(f\circ \gamma) 
\leq 
\max\Big\{0, 1-2s, \frac{1-2s}{m\kappa}\Big\}
\]
When $m<1$, the pointwise convergence \eqref{pw along curve} holds for all $f\in H^s(\mathbb R)$ if 
\[
s>
\begin{cases}
    \max\{\frac12-\frac m4,\frac{1-m\kappa}{2}\}\quad &\text{for $\kappa \in (0,\infty)$},\\
    \frac m4 \quad &\text{for $\kappa=\infty$}.
\end{cases}
\]
Furthermore, for such $s$, its Hausdorff dimension is bounded above by
\[
\dim_H\D(f\circ \gamma) 
\leq 
\begin{cases}
\max\{0,1-2s, \frac{1-2s}{m\kappa}\} \quad &\text{for $\kappa\in(0,\infty)$},\\
\max\{0,1-2s, \frac12+\frac{1-4s}{2(1-m)}\}\quad &\text{for $\kappa=\infty$}.
\end{cases}
\]
The results are sharp in the sense of the corresponding maximal inequalities. (See also Figure \ref{f:summary}.)\\

\begin{figure}[t]
\begin{center}
\begin{tikzpicture}[scale=7]


 \fill [fill=cyan!10!] (0,0)--(2/3,0)--(2/3,1/2+1/10)--(0,1/2+1/10);
 \fill [fill=magenta!10!] (2/3,0)--(1+1/5,0)--(1+1/5,1/2+1/10)--(2/3,1/2+1/10);

\draw [->] (0,0)--(1+1/5,0) node [below] at (1+1/5,0) {$m$} node [below left] at (0,0) {$O$};
\draw [->](0,0)--(0,1/2+1/10) node [left] {$s$};
\draw [dotted, name path=v] (2/3,0)--(2/3,1/2) node [below] at (2/3,0) {$1$};
\draw [Green, name path=h] (2/3,1/4)--(1+1/5,1/4);
\draw [dotted, name path=h'] (0,1/4)--(2/3,1/4);
\node [left] at (0,1/4) {$\frac14$};
\draw [dotted] (0,1/2)--(1+1/5,1/2) node [left] at (0,1/2) {$\frac12$};
\draw [blue] (0,0)--(2/3,1/4);
\draw [orange, name path=decline] (0,1/2)--(2/3,1/4);
\draw [red, name path=line] (0,1/2)--(1+1/5,1/4-1/20);
\path [name intersections={of=h and line, by={A}}];
\path [name intersections={of=v and line,by={C}}];
\draw [dotted] (A)--([yshift=-1/4*1cm]A) node [below] {$\frac{1}{2\kappa}$};

\fill [violet] (2/3,1/2) circle (1/3pt);
\fill [white] (2/3,1/4) circle (1/3pt);
\fill [white] (C) circle (1/3pt);
\draw (2/3,1/4) circle (1/3pt);
\draw (C) circle (1/3pt);
 
\end{tikzpicture}
 
 \caption{Relations between $m\in(0,\infty)$ and the critical smooth regularities $s$ that guarantee the pointwise convergence (with respect to $\d x$). The point $\textcolor{violet}{s = \frac{1}{2}}$ and the two lines $\textcolor{orange}{s = \frac{1}{2} - \frac{m}{4}}$ and $\textcolor{Green}{s = \frac{1}{4}}$ denote the least regularities required for $m \in (0,\infty)$. For each $m$, the line $\textcolor{red}{s = \frac{1 - m\kappa}{2}}$ varies and imposes additional conditions depending on $\kappa > 0$.
 Specifically, the required regularity changes at $m = \frac{1}{2\kappa}$ for $m > 1$ (if $\kappa<\frac12$) and at $\kappa = \frac{1}{2}$ for $m < 1$. When $m < 1$ and $\kappa = \infty$, the regularity is dominated by $\textcolor{blue}{s = \frac{m}{4}}$.}
 \label{f:summary}
\end{center}
\end{figure}

Lastly, we comment on another variation introduced by Sj\"ogren--Sj\"olin \cite{SS89} and Lee--Vargas with the first author \cite{CLV12}; pointwise convergence along a set of lines generated by a compact set $\Theta\subset \mathbb R^n$. When $n=1$, \cite{SS89,CLV12,Sh19} have shown that, for $m>1$, $\Theta\subset \mathbb R$ and $\gamma(x,t)=x-\theta t$ with $\theta\in \Theta$, \eqref{pw along curve} holds for all $f\in H^s(\mathbb R)$ if $s>\frac{1+\beta(\Theta)}{4}$, where $\beta(\Theta)$ denotes the Minkowski dimension of $\Theta$.  The necessity of the conditions are known in some special cases when $\Theta=\{0\},[0,1]$ whose Minkowski dimension is $0,1$, respectively. For $\Theta$ having an intermediate dimension $\beta(\Theta)\in (0,1)$ it is still open. Some higher-dimensional results are available (see \cite{Jo10,LWY20}).

In the maximal inequality perspective, however, the authors very recently proved the optimality in a sense. In fact, there exists $C>0$ such that
\begin{equation}\label{e:max fractal}
\|e^{it(-\Delta)^{\frac m2}} f(\gamma(x,t))\|_{L_x^4(0,1)L_t^\infty(0,1)L_\theta^\infty(\Theta)}
\lesssim
C
\|f\|_{H^s}
\end{equation} 
holds for all $f\in H^s(\mathbb R)$ if $s>\frac{1+\beta(\Theta)}{4}$ (due to the second author \cite{Sh19}), and this result is sharp since there exist a compact set $\Theta\subset \mathbb R$ and an initial function contained in $H^s$ such that \eqref{e:max fractal} fails when $s<\frac{1+\beta(\Theta)}{4}$. For the sharpness considering $L_x^4$ is crucial in the argument. We do not know whether it can be extended to, say, $L_x^2$ which may be more closely related to the pointwise convergence formulation.

In the case $m\in(0,1)$, there is no result in this direction as far as the authors are aware.
Modifying the argument for Theorem \ref{t:nontan}, one can show that there exists some constant $C>0$ such that \eqref{e:max fractal} holds for all $f\in H^s(\mathbb R)$ if $s>\frac12-\frac m4+\frac{m\beta(\Theta)}{4}$ and its sharpness in the sense that there exist $\Theta$ and an initial data such that \eqref{e:max fractal} fails otherwise. As a consequence, we have the following. 
\begin{theorem}\label{t:pw along lines}
Let $m\in(0,1)$ and $\Theta$ be a compact set in $\mathbb R$ whose Minkowski dimension $\beta(\Theta)\in[0,1]$. Then, the pointwise convergence along lines generated by $\Theta$ 
\[
\lim_{\substack{(y,t)\to(x,0)\\ y=x-\theta t, \ \theta\in \Theta}} e^{it(-\Delta)^{\frac m2}} f(y)
=
f(x)\quad \ae
\]
holds for all $f\in H^s(\mathbb R)$ if $s>\frac12-\frac m4+\frac{m\beta(\Theta)}{4}$. 
\end{theorem}
As we see the details in Section \ref{s:fractal}, Theorem \ref{t:pw along lines} can be generalized by replacing $\d x$ by $\d \mu$ so that one can compute an upper bound of the Hausdorff dimension of the corresponding divergence sets as
\[
\dim_H\mathfrak{D}_{s}(f\circ\gamma)
\leq
\max\left\{\frac{1-2s+m\beta(\Theta)}{m}, \frac{m\beta(\Theta)}{4s-2+m}\right\}
\]
for $s\in (\frac{2-m+m\beta(\Theta)}{4},\frac12)$, but we do not know sharpness of the estimate.\\


\subsection*{Organization}
Section \ref{s:outline} is devoted to the general framework of the proof of the sufficient conditions. We discuss the maximal estimates associated with vertical lines, non-tangential lines, and sets of lines in Section \ref{s:vertical}, \ref{s:non-tangential}, \ref{s:fractal}, respectively. Each section is divided into three parts; the proof of sufficiency, its sharpness, and technical remarks. Our main novelties are the sufficient part in Section \ref{s:vertical}, where we take a similar spirit in \cite{Sh19,CS21} but employ a non-trivial decomposition of frequency to deal with the delicate small $s$ less than $\frac14$, and the sharpness part in Section \ref{s:non-tangential}. 

\subsection*{Notation} 
Let us denote $\mathbb I$ the interval $(0,1)$.  A cut-off function over an interval $I$ is given by $\chi_{I}$ that equals to $1$ on $I$ and $0$ elsewhere.
Suppose $\psi$ is a 
Schwartz function whose Fourier transform is supported in $\{\xi\in\mathbb R:2^{-1}\leq|\xi|\leq 2\}$ that may provide a standard dyadic partition of unity, under $\psi_k(\xi)=\psi(2^{-k}\xi)$. Then, for each $k\in \mathbb N$, define $\widehat{P_kf}(\xi)=\psi_k(\xi)\widehat{f}(\xi)$.
 Regarding constants, for positive $A, B$, $A\lesssim B$ and $A\gtrsim B$ mean $A\leq CB$ and $A\geq CB$ for some positive constant $C$.  We also write $A\sim B$ for both $A\lesssim B$ and $A\gtrsim B$.

\subsection*{Acknowledgment} 
This work is supported by NRF grant no. 2022R1A4A1018904 and  RS202300239774 (Republic of Korea) (Cho), and JSPS Kakenhi grant 19H01796 and the Harmonic Analysis Incubation Research Group at Saitama University (Shiraki).
The second author would like to express his appreciation to Mitsuru Sugimoto and Neal Bez for their continuous encouragement and several inspiring discussions. He also wishes to thank Jinbong Lee for his generous hospitality at Seoul National University, where part of this research was conducted. The authors extend their gratitude to the anonymous referee for their thorough review, valuable comments, and patient efforts in improving this manuscript.



\section{Outline of proofs for sufficiency}\label{s:outline}
By invoking Frostman's lemma, the upper bound of the Hausdorff dimension follows from \eqref{e:max general}, the maximal inequality with respect to the $\alpha$-dimensional measure $\mu$ and the curve $\gamma$ given by \eqref{curve} with $\theta\in[0,\infty)$ along which the convergence is considered. For the details, the readers may consult 
 \cite{BBCR11,CL14,CS21}. 

To demonstrate \eqref{e:max general} for all $f\in H^s(\mathbb R)$ whenever $s>\frac12-s_*$ for certain $s_*\in[0,\frac12]$, we employ a standard argument presented in \cite{CLV12,Sh19,CS21}. Let $\gamma$ be a curve given by \eqref{curve} and set $T_{\gamma}f=e^{it(-\Delta)^{\frac m2}} f\circ\gamma$. By the Littlewood--Paley decomposition, \eqref{e:max general} follows from 
\[
\|T_{\gamma}f\|_{L_x^q(\mathbb I,\d\mu)L_t^\infty(\mathbb I)}
\lesssim
\lambda^{1-2s_*+\varepsilon}
\|f\|_{H^2}
\]
for $f\in L^2$ whose frequency support is contained in $\{\xi:2^{-1}\lambda\leq|\xi|\leq 2\lambda\}$ with $\lambda\ge1$ (The factor $\lambda$ represents $2^{k}$ for positive $k\in \mathbb Z$ in the decomposition). Allowing a slight abuse of notation, one may see that it is enough to show its dual form 
\begin{equation}\label{goal}
\|T_{\gamma}^*g\|_{L^2}^2
\lesssim
\lambda^{1-2s_*+\varepsilon}\|g\|_{L_x^{q'}(\d\mu)L_t^1}^2
\end{equation}
for all $g\in L_x^{q'}(\mathbb R,\d\mu)L_t^1(\mathbb R)$. Here, $T_{\gamma}^*$ is defined by 
\[
T_{\gamma}^*g(\xi)
= \psi(\lambda^{-1}\xi)
\iint\chi(x',t')e^{-i(\gamma(x',t')\xi+t'|\xi|^m)}g(x',t') \, \d\mu(x') \d t'
\]
for each $\xi\in \mathbb R$ and a fixed $\lambda\geq1$. We observe that 
\begin{align*}
&\|T_{\gamma}^*g\|_{L^2}^2\\
&=
\iint\iint
g(x,t)g(x',t')\chi(x,t)\chi(x',t')K_\lambda(\gamma(x,t)-\gamma(x',t'),t-t')\,\d\mu(x)\d\mu(x')\d t\d t',
\end{align*}
where 
$K_\lambda(x,t)=\lambda\int e^{i(\lambda x\xi+\lambda^m t|\xi|^m)}\psi^2(\xi)\,\d\xi$.
If one shows that 
\[
|K_\lambda(\gamma(x,t)-\gamma(x',t'),t-t')|\lesssim J_\lambda^\gamma(x-x')
\]
for some $ J_\lambda$ structured so that the following lemma can be applied. This yields \eqref{goal}.



\begin{lemma}[\cite{CS_RIMS}]\label{l:Y/HLS}
Let $0<\alpha\le1$, \textcolor{black}{$q\geq2$} and $\mu$ be an $\alpha$-dimensional measure. There exists a constant $C$ such that for any $b>0$,
it holds that
\begin{align}\label{e:Young special case}
&\left|\iint\iint g(x,t)h(x',t')\chi_{(0,b)}(x-x')\,\d \mu(x)\d t\d \mu(x')\d t'\right|\\
&\qquad\le Cb^{\frac{2\alpha}{q}}\|g\|_{L^{q'}_x(\d\mu)L^1_t}\|h\|_{L^{q'}_x(\d\mu)L^1_t}.\nonumber
\end{align}
Moreover, for \textcolor{black}{$0<\frac{q\rho}{2}<\alpha$} there exists a constant $C$ such that 
\begin{align}\label{e:HLS-type}
&\left|\iint\iint g(x,t)h(x',t')|x-x'|^{-\rho}\,\d \mu(x)\d t\d \mu(x')\d t'\right|\\
&\qquad\le C\|g\|_{L^{q'}_x(\d\mu)L^1_t}\|h\|_{L^{q'}_x(\d\mu)L^1_t}.\nonumber
\end{align}
Here, the integrals are taken over $(x,t,x',t')\in \mathbb I^4$.
\end{lemma}

\section{Along a vertical line}\label{s:vertical}
In this section, we prove Theorem \ref{t:concave}.
\subsection{Sufficiency}
Fix $m\in(0,1)$, $\alpha\in(0,1]$ and $q\geq2$. Our goal here is to show that \eqref{e:max} with $\gamma(x,t)=x$ holds whenever $s$ satisfies \eqref{s:vertical}.
Let us set $s_*=\min\{ \frac m4+\frac{(1-m)\alpha}{q},\frac\alpha q \}$.
By the argument in Section \ref{s:outline} it is reduced to prove
\begin{equation}\label{e:kernel}
 J_\lambda^\gamma(x)
=
\lambda
\left(
\chi_{\{|\cdot|\leq\lambda^{-\frac{qs_*}{\alpha}}\}}(x)
+
\lambda^{-2s_*+\varepsilon}|x|^{-\frac{2\alpha}{q}+\varepsilon}
\right)
\end{equation}
for $(x,t)\in \mathbb I\times \mathbb I$, arbitrary small $\varepsilon>0$, and $\lambda\geq1$. 
The first term of \eqref{e:kernel} is easily obtained from the trivial kernel estimate. Thus we may assume $|x|\geq \lambda^{-\frac{qs_*}{\alpha}}$ in the rest of the proof. In particular, note $\lambda|x|\geq1$. Denoting $\phi(\xi)=\lambda x \xi+\lambda^mt|\xi|^m$ for fixed $x,t$, we split the integral into two pieces as follows:
\[
\int e^{i\phi(\xi)}\psi(\xi)\, \d \xi
=
\int_{V_1}e^{i\phi(\xi)}\psi(\xi)\,\d\xi
+
\int_{V_2}e^{i\phi(\xi)}\psi(\xi)\,\d \xi
:=
I_1+I_2,
\]
where 
\[
V_1(x,t)
=
\{
\xi\in(2^{-1},2): 2\lambda^m |t||\xi|^{m-1}\geq \lambda^{4s_*}|x|^{\frac{4\alpha}{q}}
\},
\]
\[
V_2(x,t)
=
\{
\xi\in(2^{-1},2): 2\lambda^m |t||\xi|^{m-1}< \lambda^{4s_*}|x|^{\frac{4\alpha}{q}}
\}.
\]

One can readily deal with $I_1$. Since $|\frac{\d^2}{\d\xi^2}\phi(\xi)|\gtrsim\lambda^{4s_*}|x|^{\frac{4\alpha}{q}}$, van der Corput's lemma yields that 
\[
|I_1(x,t)|
\lesssim
(\lambda^{4s_*}|x|^{\frac{4\alpha}{q}})^{-\frac12}
\lesssim
\lambda^{-2s_*+\varepsilon}|x|^{-\frac{2\alpha}{q}+\varepsilon}.
\]
For $I_2$, by the condition from $V_2$ and $|t|<1$, we have note that
\begin{equation}\label{modify phase}
\lambda^m|t|^{\max\{1,\frac{4\alpha}{q}\}}|\xi|^{m-1}
\leq
2^{-1}\lambda^{4s_*}|x|^{\frac{4\alpha}{q}}.
\end{equation}
Also, note that 
\[
\lambda^m(\lambda^{4s_*-m}|x|^{\frac{4\alpha}{q}})^{(\max\{1,\frac{4\alpha}{q}\})^{-1}}
=
(\lambda|x|)^{\min\{\frac{4\alpha}{q},1\}}
\leq
\lambda|x|.
\]
Thus, we obtain
\begin{align*}
|\frac{\d}{\d \xi}\phi (\xi)|
\gtrsim
\lambda|x|-m\lambda^m|t||\xi|^{m-1} 
\gtrsim
\lambda|x|-\lambda^m(\lambda^{4s_*-m}|x|^{\frac{4\alpha}{q}})^{\max\{1,\frac{4\alpha}{q}\}^{-1}}
\gtrsim
\lambda|x|
\end{align*}
so that van der Corput's lemma implies that 
\[
|I_2(x,t)|
\lesssim
(\lambda|x|)^{-1}
\lesssim
\lambda^{-2s_*+\varepsilon}|x|^{-\frac{2\alpha}{q}+\varepsilon}.
\]



\qed

\subsection{Sharpness}

The proof of sharpness of \eqref{s:verline} is rather straightforward from the previous results in Section 5 of \cite{CS21,YZ21} based on the Knapp-type examples. In fact, under the usual setting $\d\mu(x)=|x|^{\alpha-1}\d x$, one may employ the initial data whose Fourier transform is given by $\lambda^{m-2}\psi(\lambda^{m-2}\xi+\lambda^m)$ for the first condition and by $\psi(\lambda^{-1}\xi)$ for the second condition.
\subsection{Remarks}
 \begin{enumerate}[(i)]
    
\item To obtain Theorem \ref{t:concave}, let $q=2$ that provides the smallest bound of $s$.

\item The proof of sufficiency is much easier when $m>1$. As the second author essentially presented in \cite{Sh19}, one can choose $V_1$ and $V_2$ trivially as
\[
V_1(x,t)
=
\{
\xi\in(2^{-1},2): 2m\lambda^m|t||\xi|^{m-1}\geq\lambda|x|
\},
\]
\[
V_2(x,t)
=
\{
\xi\in(2^{-1},2): 2m\lambda^m|t||\xi|^{m-1}<\lambda|x|
\},
\]
then applying van der Corput's lemma to each case in order to obtain $s>\max\{\frac14,\frac{1-\alpha}{2}\}$ is sufficient.

\item In the proof given above for $m\in(0,1)$, the way of division $V_1$ and $V_2$ makes the argument on $V_1$ straightforward, then $s_*$ is determined only from the one on $V_2$. In comparison with the case $m>1$, in which \eqref{e:max general} holds if $s>\max\{\frac14,\frac12-\frac{\alpha}{q}\}$ for $\alpha \in (0,1]$, $q\geq2$ and $\gamma(x,t)=x$, one can utilize the same division $V_1$ and $V_2$: Set $s_*:=\min\{\frac14,\frac{\alpha}{q}\} \leq 1$.  Instead of \eqref{modify phase}, we employ the fact that 
\[
\lambda^{4s_*}|x|^{\frac{4\alpha}{q}}
\lesssim
\lambda^{4s_*}|x|^{\frac{4\alpha}{q}\max\{1,\frac{4\alpha}{q}\}^{-1}}
=
\begin{cases}
\lambda|x|\quad & \text{if $\alpha\geq\frac q4$},\\
(\lambda|x|)^{\frac{4\alpha}{q}} \quad & \text{if $\alpha<\frac q4$}.
\end{cases}
\]
This modification is independent of $m$, and so is the whole argument.

\end{enumerate}
\section{Along a non-tangential curve/tilted line}\label{s:non-tangential}
The aim of this section is to prove the maximal inequality \eqref{e:max general}. In contrast to the previous section, there is a trade-off that sufficiency is less difficult than the one for a vertical line while the part for sharpness requires more work.

\subsection{Sufficiency}

Let us recall $q\geq2$ and set $s^*=\min\{\frac m4, \frac{m\alpha}{q}\}$. Our goal is to show 
\begin{equation}\label{e:kernelest along curve}
 J_\lambda^\gamma(x)
=
\lambda
\left(
\chi_{\{|\cdot|\leq\lambda^{-\frac{qs_*}{\alpha}}\}}(x)
+
\lambda^{-2s_*+\varepsilon}|x|^{-2s_*+\varepsilon}
\right).
\end{equation}
Again, the first term is trivial so we shall assume $|x|>\lambda^{-\frac{qs_*}{\alpha}}$ in the rest of the proof. 
The argument is somehow simpler than \cite{CS21,YZ21}.
\[
\int e^{i{\phi}(\lambda\xi)}
=
\int_{V_1}e^{i{\phi}(\lambda\xi)}\,\d\xi
+
\int_{V_2}e^{i{\phi}(\lambda\xi)}\,\d\xi
=:
I_1+I_2,
\]
where ${\phi}(\xi)=((x-x')-(t^\kappa-t'^\kappa))\xi+(t-t')|\xi|^m$, 
\[
V_1(x-x',t-t')
=
\{
\xi\in(2^{-1},2):
(\kappa+2)|t-t'|\leq |x-x'|
\}
\]
and
\[
V_2(x-x',t-t')
=
\{
\xi\in(2^{-1},2):
(\kappa+2)|t-t'|> |x-x'|
\}.
\]

Note first that 
\[
|t^\kappa-t'^\kappa|
\leq
\kappa|t-t'|
\]
for $\kappa>0$ due to the mean value theorem. For $I_1$, observe that 
\begin{align*}
|\frac{\d}{\d \xi}\phi(\xi)|
&\geq
\lambda(|x-x'|-|t^\kappa -t'^\kappa|)-m\lambda^m|t-t'||\xi|^{m-1}\\
&\geq 
\lambda|x-x'|-(\kappa+2m)\lambda|t-t'|\\
&\gtrsim
(\lambda|x-x'|)^{2s_*}.
\end{align*}
Then, one may apply van der Corput Lemma. 
For $I_2$, note
\[
|\frac{\d^2}{\d^2\xi}\phi(\xi)|
\sim
\lambda^m|t-t'|
\gtrsim
\lambda^m|x-x'|,
\]
and then apply van der Corput's lemma to complete the proof.

\subsection{Sharpness}
Let $\d\mu(x)=|x|^{\alpha-1}\d x$ as always. Then $\mu$ is $\alpha$-dimensional. Indeed, for a fixed ball $B(a,r) = (a-r,a+r)$, if $|a|<2r$ then $B(a,r) \subset B(0,3r)$. So, 
\begin{align*}
    \mu(B(a,r)) \le \int_{-3r}^{3r} |x|^{\alpha-1} dx
                =C_\alpha r^{\alpha}.
\end{align*} 
Otherwise, we may assume that $a\ge 2r$ and apply the mean value theorem to get
\begin{align*}
    \mu(B(a,r)) 
                =\frac1\alpha \Big( (a+r)^{\alpha}-(a-r)^{\alpha} \Big)
                = \frac1\alpha r (a_*)^{\alpha-1}\le C_\alpha r^{\alpha}
\end{align*} 
since $r \le a-r<a_*< a+r$ and $|x|^{\alpha-1}$ is decreasing. The sharpness of the condition $s\geq\frac12-\frac m4$ naturally comes out by considering the initial data whose Fourier transform is given by $\lambda^{m-2}\psi(\lambda^{m-2}\xi+\lambda^m)$. However, the case for the other condition requires some new ideas. For this proof let us first fix $m\in(0,1)$ and consider $f$ such that 
\[
\widehat{f}(\xi)
=
e^{i((\frac12)^\kappa\theta\xi-\frac12|\xi|^m)}\lambda^{-1}\psi(\lambda^{-1}\xi)
\]
for some large $\lambda>1$. By the change of variables; $t=\frac12+\tau$, 
\[
\sup_{t\in\mathbb I}|e^{it(-\Delta)^{\frac m2}} f(\gamma(x,t))|
=
\sup_{\tau\in [-\frac12,\frac12]}\left|\int e^{i(\lambda(x-(\tau+\frac12)^\kappa\theta)\xi+\lambda^m\tau|\xi|^m)}e^{\lambda(\frac12)^\kappa\theta\xi}\psi(\xi)\,\d\xi\right|.
\]
Taylor's expansion in $\tau$ around the origin gives
\[
(2\tau+1)^\kappa
=
\sum_{j=0}^{N-1}a_j(2\tau)^j+O(|\tau|^N),
\]
where $(a_j)_j$ are appropriate constants, and $N$ is so large that $1<mN$ holds. Now, if we set $x=2^{-\kappa}\sum_{j=1}^{N-1}a_j(2\tau)^j=:h_N(\tau)$, then it is easy to see that $h_N$ is bijection and monotone increasing on $(0,\frac{1}{100}\lambda^{-m})$. Hence, one can find a function $\tau(x)=h_N^{-1}(x)$ satisfying
\[
0=h_N^{-1}(0) \leq \tau(x) \leq h_N^{-1}(\lambda^{-m}) \lesssim \lambda^{-m}.
\]
Such pair $(x,\tau)$ leads the phase fairly small;
\[
|\lambda(x-(\tau(x)+2^{-1})^\kappa)\xi+\lambda^m|\xi|^m\tau(x))+\lambda 2^{-\kappa}\xi|
\leq2^{-1},
\]
which implies that
\[
\|e^{it(-\Delta)^{\frac m2}} f(\gamma(x,t))\|_{L^q(\mathbb I,\d\mu)L^\infty(\mathbb I)}
\gtrsim
\lambda^{-\frac{m\alpha}{q}}.
\]
Recalling $\|f\|_{H^s}\lesssim\lambda^{s-\frac12}$ and sending $\lambda$ to $\infty$, we obtain the desired conclusion.
\subsection{Remarks}
\begin{enumerate}[(i)]
\item 
If a non-tangential curve $\gamma$ is similar enough to a vertical line around $t=0$, for instance, a curve $\gamma_0(x,t)=x-e^{-\frac1t}$, one can show this curve behaves as a vertical line in this context since the term with $e^{-\frac1t}$ is negligible when $t$ is very small.

\item The same kernel estimate \eqref{e:kernelest along curve} holds for way more general curves $\gamma:\mathbb I\times \mathbb I\to \mathbb R$ that satisfy the lipschitz condition in $t$ and a certain lower bound in $x$, namely, 
\[
|\gamma(x,t)-\gamma(x,t')|\leq C_1|t-t'|,\quad t\in\mathbb I, 
\]
\[
|\gamma(x,t)-\gamma(x',t)|\geq C_2|x-x'|,\quad x\in\mathbb I
\]
for the some constants $C_1$, $C_2$ independent of $x$, $x'$, $t$, $t'$. 

\end{enumerate}

\section{Along a set of lines} \label{s:fractal}
Finally, we note a sketch proof of Theorem \ref{t:pw along lines}.
\subsection{Sufficiency}
Let $m\in (0,1)$, $q\geq2$ and $\Theta$ be a compact set of $\mathbb R$. We shall concern with the following maximal inequality in more general setting:
\begin{equation}\label{e:max fractal q}
\|e^{it(-\Delta)^{\frac m2}} f(x-\theta y)\|_{L_x^q(\mathbb I,\d\mu),L_t^\infty(\mathbb I)L_\theta^\infty(\Theta)}
\lesssim
\|f\|_{H^s}.
\end{equation}
In this case a bit of care is needed because of the extra parameter $\theta$ before using the reduction argument in Section \ref{s:outline}. First of all,  let us set 
\[
s_*
=
\min\left\{\frac m4,\frac{\alpha}{q}\right\}
\]
and recall Littlewood--Paley decomposition $f=\sum_{k\geq0}P_kf$. Then, for each $k$, decompose $\Theta$ by intervals $\Omega_{k,j}$ of length $|\Omega_{k,j}|=(2^k)^{-\frac{qs_*}{\alpha}}$. Hence, it is reduced to 
\[
\|e^{it(-\Delta)^{\frac m2}}f(x-\theta y)\|_{L_x^q(\mathbb I,\d\mu),L_t^\infty(\mathbb I)L_\theta^\infty(\Theta)}
\lesssim
\sum_{k\geq0}\biggl(\sum_{j=1}^{N_k}\|\sup_{\substack{t\in \mathbb I\\ \theta\in\Omega_{j,k}}}|e^{it(-\Delta)^{\frac m2}} P_kf(x-\theta t)|\|_{L_x^q(\mathbb I)}^q\biggr)^\frac1q,
\]
where $N_k$ denotes the smallest number of the covering.
By the argument in Section \ref{s:outline} and the fact $N_k\lesssim(2^k)^{\frac{qs_*}{\alpha}\beta+\varepsilon}$ for small $\varepsilon>0$, it is suffices to show the following.
\begin{proposition}
Let $\lambda\geq1$, $q\geq2$ and $\Omega$ be an interval of length $\lambda ^{-\frac{qs_*}{\alpha}}$. For arbitrary small $\varepsilon>0$,
\[
\|\sup_{\substack{t\in \mathbb I\\ \theta\in \Omega}}e^{it(-\Delta)^{\frac m2}} f(x-\theta y)\|_{L_x^q(\mathbb I,\d\mu)}
\lesssim
\lambda^{\frac12-s_*+\varepsilon}
\|f\|_{L^2}
\]
holds for all $f$ Fourier-supported in $\{|\xi|\sim1\}$.
\end{proposition}
The interval of $\Omega$ is chosen so in order to verify the relation
\[
|\theta-\theta'|
\leq
\lambda^{-\frac{qs_*}{\alpha}}
\lesssim
|x-x'|.
\]
We omit the further details of the proof of this proposition since it is similar to the one in the previous section containing a single line situation. The readers may also consult with \cite{Sh19,CS_RIMS}.

\subsection{Sharpness}
We show that for $m\in(0,1)$, $q\geq2$, and $\alpha=1$ there exist the initial data $f$ and the compact set $\Theta\subset \mathbb R$ such that \eqref{e:max fractal q} fails if $s<\frac12-\frac{m}{q}+\frac{m\beta(\Theta)}{q}$. Before the proof, we shall define the $r$-th Cantor sets $\C(r)$ for $r\in(0,1)$ that plays an important role. By letting $\C_0(r)=[0,1]$, the set $\C(r)$ is constructed by infinite intersection of $\C_k(r)$ that is inductively generated by removing an interval of length $r^k(1-2r)$ from the middle of each interval consisting of $\C_{k-1}(r)$. One may write $\C_k(r)=\bigcup_{j=1}^{2^k}\Omega_{k,j}$ with the interval $\Omega_{k,j}$ of length $r^k$ for each $k,j$. Note that $\C_0(r) \supset \C_1(r) \supset \dots \supset \C_k(r) \supset \dots \supset \C(r)$ and $\beta(\Theta)=-\frac{\log2}{\log r}\in (0,1)$ for $r\in (0,\frac12)$ and $\Theta=\C(r)$.

By setting $\Theta=\C(r)$ for a fixed $r\in(0,\frac12)$, the proof goes as follows: Let $\lambda_k=r^{-k}$ for each $k\in\mathbb N$ and $f$ satisfy 
\[
\widehat{f}(\xi)
=
e^{i|\xi|^m}\psi(\lambda_k^{-\frac1m}\xi).
\]
Then, by the change of variables $t\mapsto 1-\tau$
\[
\sup_{\substack{t\in\mathbb{I}\\\theta\in\Theta}}
|e^{it(-\Delta)^\frac m2} f(x+t\theta)|
=
\lambda_k^\frac1m
\sup_{\substack{\tau\in\mathbb{I}\\ \theta\in \Theta}}
\left|\int e^{i(\lambda_k^{\frac1m}(x-\theta(x)+\tau(x)\theta(x))\xi-\lambda_k\tau(x)|\xi|^m)}\psi(\xi)\,\d\xi\right|.
\]
To make the phase fairly small, choose $\theta(x)\in \Theta$ and $\tau(x)\in \mathbb{I}$ such that $|x-\theta(x)|<\lambda_k^{-1}$ for $x\in \C_k(r)\cap(\frac12,1)$ and 
\[
\tau(x)=\frac{\theta(x)-x}{\theta(x)}.
\]
Hence, it gives that
\begin{align*}
\Big(\int_{\mathfrak C_k(r)}\sup_{\substack{t\in\mathbb{I}\\\theta\in \mathfrak C(r)}}
|e^{it(-\Delta)^\frac m2} f(x+t\theta)|^q\,\d x\Big)^\frac1q
&\gtrsim
\lambda_k^\frac1m
\left(
\sum_{j=1}^{2^{k-1}}|\Omega_{k,j}|
\right)^\frac1q
\sim
2^{\frac kq}\lambda_k^{\frac1m-\frac1q}.
\end{align*}
Combining this with $\|f\|_{H^s}\lesssim \lambda_k^{\frac{s}{m}+\frac1{2m}}$, we have 
\[
\lambda_k^{\frac1m+\frac{\beta(\Theta)}{q}-\frac1q}
\lesssim
\lambda_k^{\frac sm+\frac1{2m}}
\]
since $2^k=(r^{-k})^{\beta(\theta)}=\lambda_k^{\beta(\Theta)}$. Letting $k\to\infty$ leads what we claimed.

\subsection{Remark}

In the $\alpha$-dimensional setting with $\alpha\in(0,1]$, set $\d\mu(x)=|x|^{\alpha-1}\d x$. The proof above, further applying the mean value theorem, may show $s\geq\frac12+\frac{m\alpha\beta(\Theta)}{q}-\frac{m\alpha}{q}$ is necessary for \eqref{e:max fractal q} when $m\in(0,1)$ and $\alpha\in (0,1]$.

%


\end{document}